\documentclass{siamart190516}

\usepackage{amsmath,amsfonts,amssymb}

\newtheorem{remark}{Remark}[section]

\usepackage{graphicx}                         
\usepackage{setspace}  
\usepackage{caption}
\usepackage{listings}
\usepackage{tabularx}
\usepackage{xcolor}
\usepackage{indentfirst}
\usepackage{subcaption}

\usepackage{multirow}
\usepackage{overpic}
\usepackage{diagbox}

\usepackage{verbatim}

\usepackage{enumerate}

\usepackage[normalem]{ulem}

\usepackage{tikz}

\usepackage[title]{appendix}

\usepackage{bm}

\usepackage[font=footnotesize,labelfont=bf]{caption}

\usepackage{nicefrac}

\makeatletter
\@addtoreset{equation}{section}
\makeatother

\renewcommand{\thefigure}{\arabic{section}.\arabic{figure}}
\renewcommand{\thetable}{\arabic{section}.\arabic{table}}
\renewcommand{\theequation}{\arabic{section}.\arabic{equation}}

\newcommand\dd{\mathrm{d}}
\newcommand\pp{\partial}

\newcommand\x{\bm{x}}

\newtheorem{thm}{Theorem}[section]

	\newcommand\be {\begin{equation}}
	\newcommand\ee {\end{equation}}
	\newcommand\bu {{\bf u}}

	\newcommand\dt {{\Delta t}}

\usepackage{mhchem}

\begin{document}

\title{Convergence analysis of the variational operator splitting scheme for a reaction-diffusion system with detailed balance}

\author{Chun Liu\thanks{Department of Applied Mathematics, Illinois Institute of Technology, Chicago, IL 60616, USA (\email{cliu124@iit.edu})} \and Cheng Wang\thanks{Department of Mathematics, University of Massachusetts Dartmouth, Dartmouth, MA 02747 (\email{cwang1@umassd.edu})}\and Yiwei Wang\thanks{orresponding author. Department of Applied Mathematics, Illinois Institute of Technology, Chicago, IL 60616, USA (\email{ywang487@iit.edu})} 
\and Steven M. Wise\thanks{Department of Mathematics, The University of Tennessee, Knoxville, TN 37996, USA (\email{swise1@utk.edu})}}




\maketitle 
 
\begin{abstract}
We present a detailed convergence analysis for an operator splitting scheme proposed in [C. Liu et al., \emph{J. Comput. Phys.}, 436, 110253, 2021]
for a reaction-diffusion system with detailed balance. The numerical scheme has been constructed based on a recently developed energetic variational formulation, in which the reaction part is reformulated in terms of the reaction trajectory, and both the reaction and diffusion parts dissipate the same free energy. The scheme is energy stable and positivity-preserving.
In this paper, the detailed convergence analysis and error estimate are performed for the operator splitting scheme. 
The nonlinearity in the reaction trajectory equation, as well as the implicit treatment of nonlinear and singular logarithmic terms, impose challenges in numerical analysis. To overcome these difficulties, we make use of the convex nature of the logarithmic nonlinear terms, which are treated implicitly in the chemical reaction stage. In addition, a combination of a rough error estimate and  a refined error estimate leads to a desired bound of the numerical error in the reaction stage, in the discrete maximum norm. Furthermore, a discrete maximum principle yields the evolution bound of the numerical error function at the diffusion stage. As a direct consequence, a combination of the numerical error analysis at different stages and the consistency estimate for the operator splitting results in the convergence estimate of the numerical scheme for the full reaction-diffusion system.

\end{abstract}



\section{Introduction}
Reaction-diffusion type equations have wide applications in modeling many physical and biological systems, 
such as pattern formation~\cite{hao2020spatial, kondo2010reaction, pearson1993complex}, tumor growth \cite{hawkins2012numerical, liu2018accurate, perthame2014hele}, molecular motors \cite{chipot2003variational, julicher1997modeling, wang2003robust} and active materials~\cite{prost2015active, wang2021two}, etc. 
For simplicity of presentation, we consider a reaction-diffusion system with three reactive components: 
\begin{equation}\label{Reaction-D}
  \begin{cases}
    & \pp_t a = \nabla \cdot ( D_a({\bf x}) \nabla a ) - ab + c , \\
    & \pp_t b = \nabla \cdot ( D_b({\bf x}) \nabla b ) - ab + c ,  \\
    & \pp_t c = \nabla \cdot ( D_c({\bf x}) \nabla c ) + ab - c , \\
    \end{cases} 
\end{equation}
subject to a periodic boundary condition
and the positive initial condition $$(a({\bf x}, 0), b({\bf x}, 0), c({\bf x}, 0)) = (a_0({\bf x}), b_0({\bf x}), c_0({\bf x}))\in \mathbb{R}^{3,+}.$$ Here $a$, $b$ and $c$ are concentrations of species $A$, $B$ and $C$, $D_{\alpha}({\bf x}) > 0$ ($\alpha = a, b, c$) are diffusion coefficients. The system (\ref{Reaction-D}) is associated to a chemical reaction $\ce{A + B <=>[k^+][k^-] C}$, with $k^+ = k^- = 1$.

The original reaction-diffusion system~\eqref{Reaction-D} is not a gradient flow, at least not in a direct form. As a result, standard numerical methodologies for gradient flows are not directly applicable to this system. 
Fortunately, some recent works~\cite{liero2013gradient, liu2020structure, wang2020field} have discovered that the reaction  and diffusion parts correspond to two different, but complimentary gradient flow structures. 
Although these two gradient flow structures are very different, they share exactly the same free energy. 

Indeed, let  $a^{\infty} > 0$, $b^{\infty} > 0$ and $c^{\infty} > 0$ satisfying
\begin{equation} \label{balance-1} 
k^+ a^{\infty} b^{\infty} = k^- c^{\infty},
\end{equation}
with $k^+ = k^- = 1$ in the present case, we can define the free energy $\mathcal{F}(a, b, c)$ 
as 
\begin{equation} 
\mathcal{F}(a, b, c) := \int_\Omega \left( a \Big( \ln \left( \frac{ a }{a^{\infty}} \right) - 1 \Big) 
 + b \Big( \ln \left( \frac{ b }{b^{\infty}} \right) - 1 \Big)  + c \Big( \ln \left( \frac{ c }{c^{\infty}} \right) - 1 \Big)  
 \right) \,\mathrm{d} {\bf x} .  
  \label{energy-Field} 
  \end{equation}
The corresponding chemical potentials, $\mu_A$, $\mu_B$ and $\mu_C$, for species $A$, $B$ and $C$ associated to the free energy $\mathcal{F}(a, b, c)$, can be calculated as
	\begin{equation}
  \mu_A := \frac{\delta \mathcal{F}}{\delta a} = \ln \frac{a}{a^{\infty}} , \, \, 
  \mu_B := \frac{\delta \mathcal{F}}{\delta b} = \ln \frac{b}{b^{\infty}}  , \, \, 
  \mu_C := \frac{\delta \mathcal{F}}{\delta c} = \ln \frac{c}{c^{\infty}}. 
	\end{equation}
For the reaction-only part 
\begin{equation} 
\partial_t a  = -ab + c ,  \quad \partial_t b = -ab + c ,  \quad \partial_t c = ab-c , \label{reaction-1} 
\end{equation} 
one can introduce a new variable
\[ 
R (t) = \int_0^t (ab-c) \, ds,
\] 
known as the reaction trajectory \cite{wang2020field}. 
The reaction trajectory $R(\x, t)$, which was originally introduced by de Donder~\cite{de1927affinite} as a state variable for a chemical reaction system, accounts for the number of forward reaction which has happened by time $t$.
 In turn, one gets $a = a_0 - R$, $b = b_0 - R$, $c = c_0 + R$, and the following equation of $R$ could be derived~\cite{wang2020field}: 
\begin{align} 
\ln \Big( 1 + \frac{\partial_t R}{c} \Big) = & \ln \Big( 1 + \frac{ab -c}{c} \Big)  
= \ln \Big( \frac{ab}{c} \Big)  =   \ln a + \ln b - \ln c    \nonumber 
\\
= & 
\ln (a_0 -R) + \ln (b_0 -R) - \ln (c_0 +R)  .  
\label{reaction-2} 
\end{align} 
The free energy can be written in terms of $R$, specifically, $\mathcal{F}(a, b, c) = \tilde{\mathcal{F}}(R)$, and it is easy to see that
\begin{equation}\label{Equation_R}
\ln \Big( 1 + \frac{\partial_t R}{c} \Big) =  \ln (a_0 -R) + \ln (b_0 -R) - \ln (c_0 +R)  = - \frac{\delta\tilde{\mathcal{F}}}{\delta R}	.
\end{equation}
Therefore, the following energy dissipation law is available: 
\begin{align*}
\frac{\dd}{\dd t} \tilde{\mathcal{F}} (R) = \left(\partial_t R, \frac{\delta\tilde{\mathcal{F}}}{\delta R}\right) = - \left( c \frac{\partial_t R}{c}, \ln \Big( 1 + \frac{\partial_t R}{c} \Big) \right) \le 0,
\end{align*}
provided that $\frac{\partial_t R}{c} > -1$. In other words, the reaction part becomes a generalized gradient flow in terms of $R$, which is decidedly different from the standard $L^2$ or $H^{-1}$ gradient flow structures. The monotonicity of $\ln ( 1 + \frac{\partial_t R}{c})$ (in terms of $\partial_t R$) for $c > 0$ plays an important role in the dissipation mechanism. 

In the meantime, it is observed that the diffusion-only part 
\begin{equation} 
  \pp_t a = \nabla \cdot ( D_a({\bf x}) \nabla a ) , \quad 
  \pp_t b = \nabla \cdot ( D_b({\bf x}) \nabla b ) ,  \quad 
  \pp_t c = \nabla \cdot ( D_c({\bf x}) \nabla c ) ,  \label{diffusion-1} 
\end{equation}  
could be rewritten as the following $H^{-1}$ gradient flow, with non-constant mobility: 
\begin{equation} 
  \pp_t a = \nabla \cdot ( D_a({\bf x}) a \nabla \mu_A ) , \quad 
  \pp_t b = \nabla \cdot ( D_b({\bf x}) b \nabla \mu_B ) ,  \quad 
  \pp_t c = \nabla \cdot ( D_c({\bf x}) c \nabla \mu_C ) .  \label{diffusion-2} 
\end{equation} 
This gradient flow structure is similar to that of the Poisson-Nernst-Planck (PNP) system~\cite{liu2020positivity, QianWangZhou_JCP21}. 
As a consequence, the overall system satisfies the energy-dissipation law \cite{wang2020field}:
\begin{equation}\label{ED1}
  \begin{aligned}
\frac{\dd}{\dd t} \mathcal{F}(a, b, c)  = - \int_\Omega  & D_a({\bf x}) a |\nabla \mu_A|^2 + D_b({\bf x}) b |\nabla \mu_B|^2 + D_c({\bf x}) c |\nabla \mu_C|^2  \\
& + \pp_t R \ln \Big( 1 + \frac{\partial_t R}{c} \Big)   \dd {\bf x} \leq  0. \\
  \end{aligned}
\end{equation} 

\begin{remark}
  There have been many existing works aiming to establish a variational structure of reaction-diffusion systems \cite{anderson2015lyapunov, desvillettes2006, desvillettes2017trend, glitzky2013gradient, haskovec2018decay, mielke2011gradient, mielke2013thermomechanical, liero2013gradient, shear1967analog,wei1962axiomatic, wang2020field}. The condition~\eqref{balance-1} is known as the \emph{detialed balance} condition, which guarantees the existence of the free energy \cite{anderson2015lyapunov, desvillettes2017trend}.  
  We call (\ref{ED1}) as the energetic variational formulation for the reaction-dfifusion system, which can be used to model the coupling between a general reversible reaction newtwork and other mechanical process, such as general diffusions. We refer the interested reader to \cite{liu2020structure, wang2020field} for the energetic variational formulation for more general cases.
\end{remark}

 

%

Although the reaction and diffusion parts in (\ref{Reaction-D}) correspond to very different gradient flow structures, their free energy functionals are exactly the same. If one applies a standard numerical discretization to the original reaction-diffusion system~\eqref{Reaction-D}, the variational structure may be lost. In turn, either the theoretical property or the logarithmic energy stability could not be justified. This fact motivates the development of the operator splitting scheme \cite{liu2020structure}, in which the reaction stage is solved in terms of the reaction trajectory $R$ and  both stages dissipate the same discrete energy. The unique solvability, positivity-preserving property and energy stability have been theoretically established for the proposed operator splitting scheme. However, its convergence analysis and error estimate remain open, in which the primary difficulty comes from the nonlinear and singular nature of the logarithmic terms. The implicit treatment of these nonlinear and singular logarithmic terms are crucial to enforce the positivity of the numerical solution, as well as the energy stability analysis, while it has posed a great challenge in the theoretical justification of the convergence analysis. Also see the related works~\cite{chen19b, dong20b, dong19a, dong20a, Yuan2021a} for the Cahn-Hilliard equation with Flory-Huggins energy potential, as well as~\cite{liu2020positivity, QianWangZhou_JCP21} for the Poisson-Nernst-Planck system, \cite{duan20a} for the  porous medium equation, \cite{ZhangJ2021} for a liquid film droplet model, etc.

In this paper we provide a detailed convergence analysis and error estimate for the operator splitting scheme, proposed in \cite{liu2020structure} and applied to the reaction-diffusion system~\eqref{Reaction-D}. A careful consistency estimate for the splitting process, as well as the temporal discretization at each stage, gives an $O ( \dt)$  truncation error. In addition, the centered difference spatial discretization at the diffusion stage implies an $O (h^2)$ truncation error. To overcome the subtle difficulty associated with the singularity, we make use of the convex nature of the logarithmic nonlinear terms, which are implicitly treated in the chemical reaction stage. In addition, a combination of rough error estimate and refined error estimate are performed in the reaction stage, which in turn leads to a desired bound of the numerical error in the discrete maximum norm. Moreover, a careful application of discrete maximum principle yields the evolution bound of the numerical error function at the diffusion stage. Therefore, a combination of the numerical error analysis at different stages results in the convergence estimate of the numerical scheme for the full reaction-diffusion system, in the discrete maximum norm. 

The rest of this paper is organized as follows. 
The positive-preserving and energy stable operator splitting scheme for the reaction-diffusion system (\ref{Reaction-D}) is reviewed in Section~\ref{sec: numerical scheme}. The optimal rate convergence analysis and error estimate are presented in Section~\ref{sec: convergence}. A numerical result is given in Section 4, which validate the theoretical analysis. Finally, some concluding remarks are made in Section~\ref{sec: conclusion}.

 \section{Review of the operator splitting numerical scheme}  \label{sec: numerical scheme} 
In this section, we give a brief review to the operator splitting method proposed in \cite{liu2020structure}, which is based on the energetic variational formation (\ref{ED1}) of equation \eqref{Reaction-D}. 
Let $\bu({\bf x}, t)= (a({\bf x}, t), b({\bf x}, t), c({\bf x}, t))^{\rm T}$, the reaction-diffusion system (\ref{Reaction-D}) can be represented as
 \begin{equation}
 \bu (\x, t) = \mathcal{A} \bu  + \mathcal{B} \bu , 
 \end{equation}
 where  $\mathcal{A}$ and $\mathcal{B}$ are reaction operator and diffusion operator, respectively.
 As mentioned earlier, the key point of designing an energy-stable, positivity-preserving numerical scheme for the reaction part is to discretize the reaction trajectory equation (\ref{Equation_R}) directly. 
 We present the numerical algorithm on the computational domain $\Omega = (0,1)^3$ with  periodic boundary conditions and use a finite difference method as a spatial discretization. The spatial mesh size is set as $\Delta x = \Delta y = \Delta z = h = \frac{1}{N}$, where $N$ is the spatial mesh resolution throughout this paper. In particular, $f_{i,j,k}$ stands for the numerical value of $f$ at the cell centered mesh points $( ( i + \frac12 ) h, (j + \frac12 ) h, ( k + \frac12 ) h )$, so that the discrete summation could be easily defined over $\Omega$. The discrete gradient, divergence and Laplacian operators, given by $\nabla_h$, $\nabla_h \cdot$ and $\Delta_h$, are defined based on the standard centered difference approximation. The discrete $L^2$ inner product between two grid functions $f$ and $g$, as well as the discrete $L^2$ norm, are defined as 
 \begin{equation} 
    \langle f , g \rangle := h^3 \sum_{i,j,k=1}^{N} f_{i,j,k} g_{i,j,k} ,  \quad  
    \| f \|_2 := ( \langle f , f \rangle )^\frac12 . 
    \label{defi-inner product-1} 
 \end{equation}  
 As an application, the discrete energy of a numerical solution $(a,b,c)$ is introduced as 
 \begin{equation} 
   \mathcal{F}_h (a, b , c): = \Big\langle a \Big( \ln \left( \frac{ a }{a^{\infty}} \right) - 1 \Big) 
    + b \Big( \ln \left( \frac{ b }{b^{\infty}} \right) - 1 \Big)  + c \Big( \ln \left( \frac{ c }{c^{\infty}} \right) - 1 \Big)  , {\bf 1} \Big\rangle . \label{defi-discrete energy-1} 
 \end{equation} 
 In addition to the discrete $\| \cdot \|_2$ norm, the discrete maximum norm is defined as follows: 
 \begin{equation} 
   \| f \|_\infty := \max_{1\le i,j,k\le N }\left| f_{i,j,k}\right| . 
   \label{defi-maximum norm-1} 
 \end{equation}

Based on the energy-dissipation law~\eqref{ED1}, the operator splitting scheme for equation (\ref{Reaction-D}) can be formulated as follows:
Given $a^n$, $b^n$, $c^n$, with $a^n, b^n, c^n > 0$ at each mesh point. We update $a^{n+1}$, $b^{n+1}$, $c^{n+1}$, via the following two stages. 

\noindent 
{\bf Stage 1.} \, First, we set $R^n \equiv 0$, and solve
	\begin{equation} 
\ln \left( \frac{R^{n+1} - R^n}{ c^n \dt} + 1 \right) 
 = \ln \left( \frac{a^n - R^{n+1}}{a^{\infty}} \right) 
 + \ln \left( \frac{b^n - R^{n+1}}{b^{\infty}} \right) - \ln \left( \frac{c^n + R^{n+1}}{c^{\infty}} \right) ,  
 	\label{splitting-Field-A-1} 
	\end{equation} 
  at each mesh point.
By a careful analysis based on the convexity of the logarithmic function, one can show that there exists a unique solution $R^{n+1}$ such that $a^n - R^{n+1} > 0$, $b^n - R^{n+1} > 0$, $c^n + R^{n+1} > 0$, and $R^{n+1} - R^n + c^n \dt > 0$. In turn, we denote 
\begin{equation} 
   a^{n+1,*} := a^n - R^{n+1} , \, \, \, 
   b^{n+1,*} := b^n - R^{n+1} , \, \, \, 
   c^{n+1,*} := c^n + R^{n+1} .  \label{splitting-Field-A-3}   
\end{equation}    
Furthermore, the following energy dissipation property has been established \cite{liu2020structure}:
	\begin{equation} 
  \mathcal{F}_h (a^{n+1,*}, b^{n+1,*}, c^{n+1,*}) \le \mathcal{F}_h (a^n, b^n, c^n ) . 
     \label{splitting-Field-A-4}     
	\end{equation} 

\noindent 
{\bf Stage 2.} \,The intermediate variables $a^{n+1,*}$, $b^{n+1,*}$, $c^{n+1,*}$ have been proved to be positive at each mesh point. Next, we update $a^{n+1}$, $b^{n+1}$ and $c^{n+1}$ by the standard implicit Euler scheme
	\begin{equation}
	\label{splitting-Field-B-1}
	\left\{
  	\begin{split}
 \frac{a^{n+1} - a^{n+1, *}}{\dt} &=  \nabla_h \cdot ( D_a  \nabla_h a^{n+1}), 
	\\ 
 \frac{b^{n+1} - b^{n+1, *}}{\dt} &=  \nabla_h \cdot ( D_b  \nabla_h b^{n+1}), 
	\\ 
 \frac{c^{n+1} - c^{n+1, *}}{\dt} &=  \nabla_h \cdot ( D_c  \nabla_h c^{n+1}), 
	\\ 
  	\end{split}
  	\right.
	\end{equation}
where $\nabla_h$ and $\nabla_h \cdot$ are discrete gradient and divergence operators. The positivity and energy stability of the implicit Eulerian scheme has been proved in \cite{liu2020structure}, i.e.,
\begin{eqnarray} 
  &&  
  a^{n+1} , \, b^{n+1} , \, c^{n+1} > 0   \quad 
  \mbox{(point-wise)}  ,  \label{splitting-Field-B-5}   
\\
  && 
    \mathcal{F}_h (a^{n+1}, b^{n+1}, c^{n+1})  
    \le \mathcal{F}_h (a^{n+1,*}, b^{n+1,*}, c^{n+1,*})  . 
     \label{splitting-Field-B-6}     
\end{eqnarray} 
A combination of~\eqref{splitting-Field-A-4} and \eqref{splitting-Field-B-6} results in 
\begin{eqnarray} 
    \mathcal{F}_h (a^{n+1}, b^{n+1}, c^{n+1})  \le \mathcal{F}_h (a^n, b^n, c^n)  . 
     \label{splitting-Field-B-7}     
\end{eqnarray} 
Therefore, we arrive at the following theoretical result for the operator splitting scheme. 
  
\begin{thm}[\cite{liu2020structure}] \label{Field-operator splitting} 
  Given $a^n$, $b^n$, $c^n$, with $a^n_{i,j,k}, b^n_{i,j,k}, c^n_{i,j,k} > 0$, $\forall \, 1 \le i, j, k \le N$, there exists a unique solution $a^{n+1}$, $b^{n+1}$, $c^{n+1}$, with discrete periodic or Neumann boundary conditions, for the operator splitting numerical scheme~(\eqref{splitting-Field-A-1} combined with \eqref{splitting-Field-B-1}). The point-wise positivity is ensured: $0 <  a^{n+1}_{i,j,k}, b^{n+1}_{i,j,k}, c^{n+1}_{i,j,k}$, $\forall \, 1 \le i, j, k \le N$. In addition, we have the energy dissipation estimate: $\mathcal{F}_h (a^{n+1}, b^{n+1}, c^{n+1}) \le \mathcal{F}_h (a^n, b^n, c^n)$, so that ${\mathcal F}_h (a^n, b^n, c^n) \le {\mathcal F}_h (a^0, b^0, c^0)$.  
\end{thm}

\section{Optimal rate convergence analysis and error estimate}  \label{sec: convergence}

Numerical results in \cite{liu2020structure} indicate that the operator splitting scheme can achieve first-order accuracy in time and second-order accuracy in space. However, a theoretical justification of the convergence analysis turns out to be a challenging subject, due to the nonlinear and singular nature in the reaction part. The main theoretical result of this paper is the following convergence theorem. 

\begin{thm}
	\label{thm:convergence}
Given positive initial data $a_0, b_0, c_0  \in C^4_{\rm per}(\Omega)$, suppose the exact solution for the reaction-diffusion system~\eqref{Reaction-D}, denoted $(a_e, b_e, c_e)$, is of regularity class $[\mathcal{R}]^3$, where
	\begin{equation}
\mathcal{R} := C^2 \left(0,T; C_{\rm per}(\Omega)\right) \cap C^1 \left(0,T; C^1_{\rm per}(\Omega)\right) \cap L^\infty \left(0,T; C^4_{\rm per}(\Omega)\right).
	\label{assumption:regularity.1}
  \end{equation}
Then, provided $\dt$ and $h$ are sufficiently small,  we have
	\begin{equation}
\| a_e^n - a^n \|_\infty + \| b_e^n - b^n \|_\infty + \| c_e^n - c^n \|_\infty  
  \le C ( \dt + h^2 ) , \label{convergence-0}
	\end{equation}
for all positive integers $n$, such that $t_n=n\dt \le T$, where $C>0$ is independent of $\dt$ and $h$, $a_e^n$, $b_e^n$ and $c_e^n$ are exact solutions at $t^n$.
  \end{thm}
  

As a consequence of the regularity assumption~\eqref{assumption:regularity.1}, the following bound is available for the exact solution: 
\begin{equation} 
  \| \bu_e \|_{C^2 (0,T; C^0)} \le C_0 , \quad 
  \| \bu_e ( \cdot , t) \|_{C^4 (\Omega)} \le C_0 ,  \, \, \, \forall t \ge 0 .  
  \label{assumption:exact bound}
\end{equation}  
In particular, there exists a constant $C_0$ such that
\begin{equation}
  \sup_{{\bf x}, t} \max \{ | a_e ({\bf x}, t) | + | \partial_t a_e ({\bf x}, t) | , 
   | b_e ({\bf x}, t) | + | \partial_t b_e ({\bf x}, t) | ,  
   | c_e ({\bf x}, t) | + | \partial_t c_e ({\bf x}, t) |  \} \leq C_0 . 
\end{equation}

In addition, the following separation property is also assumed for the exact solutions: 
\begin{equation} 
   a_e ({\bf x}, t) , \, b_e ({\bf x}, t) , \, c_e ({\bf x}, t)  \ge \epsilon_0 ,  \quad  \exists \, \epsilon_0 > 0 .  
	\label{assumption:separation}    
\end{equation} 
In fact, this assumption is necessary to ensure the regularity requirement~\eqref{assumption:regularity.1} for the exact solutions, because of the $\ln \frac{a}{a^{\infty}}$, $\ln \frac{b}{b^{\infty}}$, $\ln \frac{c}{c^{\infty}}$ terms appearing in the free energy. In fact, such a separation property has already been established for the 2-D Cahn-Hilliard equation with Flory-Huggins energy potential~\cite{abels07, debussche95, elliott96b, Giorgini17a, miranville04}, and this property is expected to hold for the reaction-diffusion system~\eqref{Reaction-D} in the energetic variational formulation.

\subsection{Consistency analysis for the operator splitting scheme} 
We first perform a consistency analysis \cite{arnold2015stability} for the operator splitting scheme.
Given $\bu_e^n = (a_e^n, b_e^n, c_e^n)^T$, with the regularity assumption~\eqref{assumption:regularity.1} and separation assumption~\eqref{assumption:separation} satisfied, we introduce $\bu_e^{n+1,*} = (a_e^{n+1,*} , b_e^{n+1,*} , c_e^{n+1,*} )^T$ as the exact update of the first stage equation: $\partial_t \bu = \mathcal{A} \bu$, over the time interval ($t^n, t^{n+1})$, with initial data $\bu_e^n$. In other words, $\bu_e^{n+1,*} = (a_e^{n+1,*} , b_e^{n+1,*} , c_e^{n+1,*} )^T$ is the exact solution at $t = t^{n+1}$ for the reaction-only equation
\begin{equation}\label{Reaction_eq_n}
  \begin{cases}
   & \partial_t \bu = \mathcal{A} \bu, \\
   & \bu ( {\bf x} , t^n) = \bu_e^n ({\bf x}) .  
  \end{cases}
\end{equation} 
Meanwhile, as mentioned in the previous section, equation \eqref{Reaction_eq_n} can be reformulated as an equation of the reaction coordinate $R_e(\x, t)$ over the time interval $(t^n, t^{n+1})$, since
\begin{eqnarray*} 
  a_e ({\bf x}, t) = a_e^n ({\bf x}) - R_e ({\bf x}, t)  , \, \, 
  b_e ({\bf x}, t) = b_e^n (\x) - R_e ({\bf x}, t)  , \, \, 
  c_e ({\bf x}, t) = c_e^n (\x) + R_e ({\bf x}, t)  .  \label{consistency-2} 
\end{eqnarray*}   
The equation for $R_e$ is given by
\begin{equation}\label{eq_Re}
  \begin{cases}
  & \ln \Big( \frac{\partial_t R_e}{ c_e^n +R_e } + 1 \Big) 
  = \ln ( \frac{a_e^n - R_e}{a^{\infty}} ) + \ln ( \frac{b_e^n - R_e}{b^{\infty}} ) - \ln ( \frac{c_e^n + R_e}{c^{\infty}} ) ,  \\
  & R_e ( \cdot , t^n) \equiv 0. \\    
  \end{cases}
  \end{equation} 
 Moreover, we have
\begin{eqnarray} 
   a_e^{n+1,*}  =  a_e^n - R_e^{n+1} , \quad 
   b_e^{n+1,*}  =  b_e^n - R_e^{n+1} , \quad 
   c_e^{n+1,*}  =  c_e^n + R_e^{n+1} .   \label{consistency-4}          
\end{eqnarray}  

By a careful Taylor expansion in time, one can show that the exact equation~\eqref{eq_Re} can be approximated by the temporal discretization~\eqref{splitting-Field-A-1} with $O(\Delta t)$ accuracy:
	\begin{equation} 
\ln \Big( \frac{R_e^{n+1} - R_e^n}{ c_e^n + R_e^n } + 1 \Big) = \ln \left( \frac{a_e^n - R_e^{n+1}}{a^{\infty}} \right) 
 + \ln \left( \frac{b_e^n - R_e^{n+1}}{b^{\infty}} \right) - \ln \left( \frac{c_e^n + R_e^{n+1}}{c^{\infty}} \right) + \tau_0^{n+1} ,  \label{consistency-4-1} 
	\end{equation}  
where $R_e^n \equiv 0$ and $| \tau_0^{n+1} | \le C \dt$ is the local truncation error, at a point-wise level.
 The consistency estimate~\eqref{consistency-4-1} could be rewritten as the following equation after an exponential transform: 
\begin{eqnarray} 
  \frac{R_e^{n+1} - R_e^n}{c_e^n + R_e^n}
    =    \frac{(a_e^n - R_e^{n+1}) (b_e^n - R_e^{n+1})}{c_e^n +R_e^{n+1}} - 1 
    + \tau_1^{n+1} ,
     \label{consistency-11} 
\end{eqnarray}  
where $| \tau_1^{n+1} | \le C \dt$ due to the fact that ${\rm e}^{\tau_0^{n+1}} - 1 = O (\dt)$ for $\tau_0^{n+1} = O (\dt)$.

In the second stage, given $\bu_e^{n+1,*} = (a_e^{n+1,*} , b_e^{n+1,*} , c_e^{n+1,*} )^T$, we denote the exact update of by $\bu_e^{n+1,**} = (a_e^{n+1,**} , b_e^{n+1,**} , c_e^{n+1,**} )^T$, i.e., $\bu_e^{n+1,**}$ is the exact solution for the linear diffusion equation
	\begin{equation} 
  	\begin{aligned}
&  \partial_t \bu = \mathcal{B} \bu = \nabla \cdot ( \mathcal{D} (\x) \nabla \bu ), 
	\\  
& \bu ( \cdot , t^n) = \bu_e^{n+1,*},  
  	\end{aligned}
   	\label{consistency-5} 
	\end{equation} 
at $t = t^{n + 1}$. By a careful Taylor expansion associated with the operator splitting $\bu_e^{n+1,**} = {\rm e}^{\mathcal{B} \dt} {\rm e}^{\mathcal{A} \dt} \bu_e^n$, one can show that
\begin{equation} 
   \bu_e^{n+1,**} - \bu_e^{n+1} = O (\dt^2) . 
    \label{consistency-6} 
\end{equation}    
On the other hand, an application of implicit Euler temporal discretization to the diffusion equation system~\eqref{consistency-5} implies the following consistency estimate 
\begin{equation} 
  \frac{\bu_e^{n+1,**} - \bu_e^{n+1,*} }{\dt} 
  = \nabla \cdot ( \mathcal{D} (\x) \nabla  \bu_e^{n+1,**} ) + \tau_{2,t}^{n+1,(1)} ,  
      \label{consistency-7} 
\end{equation}   
where $| \tau_{2,t}^{n+1,(1)} | \le C \dt$ is the local truncation error.
In turn, its combination with~\eqref{consistency-6} yields 
\begin{equation} 
  \frac{\bu_e^{n+1} - \bu_e^{n+1,*} }{\dt} 
  = \nabla \cdot ( \mathcal{D} (\x) \nabla  \bu_e^{n+1} )  + \tau_{2,t}^{n+1} ,    
      \label{consistency-8} 
\end{equation}  
where $| \tau_{2,t}^{n+1} | \le C \dt$. Furthermore, the centered difference approximation for $\bu_e^{n+1}$ leads to the following truncation error estimate: 
\begin{equation} 
  |  \nabla \cdot ( \mathcal{D} (\x) \nabla  \bu_e^{n+1} ) 
  - \nabla_h \cdot ( \mathcal{D}  \nabla_h  \bu_e^{n+1} ) | \le C h^2 ,  \quad 
  \mbox{point-wise on the mesh} . 
   \label{consistency-9} 
\end{equation}  
Then we obtain the consistency estimate for the second stage: 
\begin{equation} 
  \frac{\bu_e^{n+1} - \bu_e^{n+1,*} }{\dt} = 
  \nabla_h \cdot ( \mathcal{D}  \nabla_h  \bu_e^{n+1} ) + \tau_2^{n+1} ,   
      \label{consistency-10} 
\end{equation}  
where $| \tau_2^{n+1} | \le C ( \dt + h^2 )$.

In summary, we have the consistency analysis for the operator splitting scheme
	\begin{eqnarray} 
	&& 
\frac{R_e^{n+1} - R_e^n}{(c_e^n + R_e^n) \dt} 
    = \frac{(a_e^n - R_e^{n+1}) (b_e^n - R_e^{n+1})}{c_e^n +R_e^{n+1}} - 1
    + \tau_1^{n+1},  \quad  R_e^n = 0, 
	\label{consistency-12-1} 
	\\
&&   a_e^{n+1,*} =  a_e^n - R_e^{n+1} ,  \quad  
     b_e^{n+1,*}  =  b_e^n - R_e^{n+1} , \quad  
    c_e^{n+1,*}  =  c_e^n + R_e^{n+1} ,    
	\label{consistency-12-3}    
	\\
  &&
   \frac{\bu_e^{n+1} - \bu_e^{n+1,*} }{\dt} = 
   \nabla_h \cdot (  \mathcal{D}  \nabla_h  \bu_e^{n+1} )  + \tau_2^{n+1} ,  
      \label{consistency-12-4} 
  \end{eqnarray}  
where
	\begin{equation}
| \tau_1^{n+1} | \le C \dt \quad \mbox{and} \quad | \tau_2^{n+1} | \le C ( \dt + h^2 ).
	\end{equation}
Of course, the local truncation error is of order $O(\dt + h^2)$. 
	
\subsection{Error estimate in the first stage} 
We first perform the error estimate in the reaction stage. Define the point-wise error functions:
	\begin{equation} 
\begin{aligned} 
  & 
  e_a^k := a_e^k - a^k , \, \, \, e_b^k := b_e^k - b^k , \, \, \, 
  e_c^k := c_e^k - c^k , \,  \, \, e_R^k := R_e^k - R^k , 
\\
  &  
  e_a^{n+1,*} := a_e^{n+1,*} - a^{n+1,*} , \, \, \, 
  e_b^{n+1,*} := b_e^{n+1,*} - b^{n+1,*} , \, \, \, 
  e_c^{n+1,*} := c_e^{n+1,*} - c^{n+1,*} ,    
\end{aligned}   \label{error function-2}
\end{equation}
for any $k \ge 0$, $n \ge 0$. The numerical scheme for the reaction stage~\eqref{splitting-Field-A-1} may, of course, be rewritten by an exponential transform  
	\begin{equation} 
\frac{R^{n+1} - R^n}{(c^n + R^n) \dt} = \frac{(a^n - R^{n+1}) (b^n - R^{n+1})}{c^n +R^{n+1}} -1 ,  \quad R^n \equiv 0 . 
	\label{splitting-Field-A-rewritten} 
	\end{equation}  

Subtracting the rewritten scheme~\eqref{splitting-Field-A-rewritten} from the consistency estimate~\eqref{consistency-12-1} and rearranging terms  yields
	\begin{equation}
\frac{e_R^{n+1}}{c^n \dt}  =  q_0^{n+1} e_c^n - (q_1^{n+1} + q_2^{n+1} + q_3^{n+1} ) e_R^{n+1} +  q_1^{n+1} e_a^n + q_2^{n+1} e_b^n - q_3^{n+1}  e_c^n  +\tau_1^{n+1},
	\label{error equation-A-1}    
	\end{equation} 
where
	\begin{equation}
	\label{error equation-A-2}
	\begin{split}
& q_0^{n+1} := \frac{R_e^{n+1}}{ c_e^n \cdot c^n \dt}, \quad q_1^{n+1} := \frac{b^n - R^{n+1}}{c^n +R^{n+1}}  ,  \quad   
	\\
& q_2^{n+1} := \frac{a_e^n - R_e^{n+1}}{c^n +R^{n+1}}  , \quad q_3^{n+1} := \frac{(a_e^n - R_e^{n+1}) (b_e^n - R_e^{n+1})}{(c^n + R^{n+1} ) (c_e^n +R_e^{n+1})}  .    
 	\end{split} 
	\end{equation}

	\begin{remark}
We observe that
	\begin{equation}
    \begin{aligned}
\frac{AB}{C} - \frac{(A + \xi_A) (B + \xi_B) }{ C + \xi_B }  & = \frac{AB \xi_C}{ (C + \xi_C) C} - \frac{\xi_A \xi_B}{C + \xi_C} - \frac{A \xi_B + B \xi_A}{C + \xi_C}  
	\\
& = - \frac{B + \xi_B}{C + \xi_C} \xi_A - \frac{A}{C + \xi_C} \xi_B + \frac{AB}{C(C + \xi_C)} \xi_C . 
	\end{aligned}
	\end{equation}
By taking $A = a_e^n - R_e^{n+1}$, $B = b_e^n - R_e^{n+1}$, $C = c_e^n + R_e^{n+1}$, $\xi_A = - e_a^n + e_R^{n+1}$, $\xi_B = - e_b^n + e_R^{n+1}$ and $\xi_C = - e_c^n - e_R^{n+1}$, we can obtain \eqref{error equation-A-1}.
	\end{remark}

The error evolutionary equation~\eqref{error equation-A-1} could be rewritten as 
	\begin{equation} 
M^{n+1}  e_R^{n+1} =   q_0^{n+1} e_c^n +(q_1^{n+1} e_a^n + q_2^{n+1} e_b^n - q_3^{n+1}  e_c^n ) + \tau_1^{n+1} ,   
	\label{convergence-A-3-1}  
	\end{equation} 
where $M^{n+1}$ is defined by
	\begin{equation} 
	\label{convergence-A-3-2} 
M^{n+1} := \frac{1}{c^n \dt} +  (q_1^{n+1} + q_2^{n+1} + q_3^{n+1} )  .
	\end{equation} 
To proceed with the nonlinear analysis, we first make the following \emph{a-priori} assumption for the previous time step: 
	\begin{equation} 
\| e_a^n \|_\infty  \le \dt^\frac12  + h , \quad   \| e_b^n \|_\infty  \le \dt^\frac12 + h , \quad \| e_c^n \|_\infty  \le \dt^\frac12 + h .  \label{a priori-1} 
	\end{equation}   
Such an a-priori assumption will be recovered by the optimal rate convergence analysis at the next time step, as demonstrated later.
 
A direct consequence of the assumption (\ref{a priori-1}) gives the following bound and separation property for the numerical solution at the previous time step: 
	\begin{equation} 
	\begin{aligned}  
& | a^n | \le | a_e^n | + | e_a^n | \le C_0 + 1 := C_1 ,  \quad | b^n | \le C_1 ,  \quad   | c^n | \le C_1  , 
	\\ 
& a^n  \ge a_e^n  - | e_a^n |  \ge \frac{\epsilon_0}{2} ,  \quad  b^n \ge \frac{\epsilon_0}{2} ,  \quad   c^n \ge \frac{\epsilon_0}{2} ,    
	\end{aligned}  
	\label{a priori-2}   
	\end{equation}        
provided that 
	\[
\dt^\frac12 , \, h \le \min \left( \frac{\epsilon_0}{4}, \frac12 \right).
	\] 
Here we have made use of the functional bound~\eqref{assumption:exact bound} and the separation property~\eqref{assumption:separation} for the exact solution. 

Due the positivity-preserving property for both the exact solution and the numerical solution (coming from Theorem~\ref{Field-operator splitting}), it is obvious that 
	\[ 
a_e^n - R_e^{n+1} > 0 , \quad  b^n - R^{n+1} > 0 , \quad  b_e^n - R_e^{n+1} > 0 , \quad  c^n + R^{n+1} > 0 , \quad c_e^n + R_e^{n+1} > 0 , 
	\]
which, in turn, implies that
	\begin{equation}
q_1^{n+1} > 0 , \quad  q_2^{n+1} > 0 , \quad  q_3^{n+1} > 0 . 
  	\label{convergence-A-1} 
	\end{equation} 
Meanwhile, the $C^2([0,T])$ bound for the exact solution $R_e$ indicates that $| \frac{R_e^{n+1}}{\dt} | \le C_0$. The separation estimates for the exact and numerical solutions, given by~\eqref{assumption:separation} and \eqref{a priori-2}, respectively, lead to $0 < \frac{1}{c^n \cdot c_e^n}  < \frac{4}{\epsilon_0^2}$. 
In turn, $q_0^{n+1}$ is uniformly bounded by
	\begin{equation} 
| q_0^{n+1} | \le \frac{4 C_0}{\epsilon_0^2} := C_2 .  
	\label{convergence-A-2-3} 
	\end{equation} 

A rough error estimate on $\| e_R^{n+1} \|_{\infty} \lesssim \Delta t^{1/2} + h $ can be obtained based on the following simple estimates: 
	\begin{eqnarray} 
&& M^{n+1} \ge \frac{1}{(c^n + R^n) \dt}   \Longrightarrow 0 < \frac{1}{M^{n+1}} \le  c^n \dt \le C_1 \dt ,  
	\label{convergence-A-4-1}     
	\\
  && 
M^{n+1} \ge  q_1^{n+1} + q_2^{n+1} + q_3^{n+1}  \Longrightarrow  0 < \frac{ q_1^{n+1} + q_2^{n+1} + q_3^{n+1}}{M^{n+1}} 
   \le 1 ,  
	\label{convergence-A-4-2}   
	\\
  &&
 \Big|  \frac{q_0^{n+1}}{M^{n+1}} \Big|   \le 
\Big | \frac{C_2}{\frac{1}{ c^n \dt} }  \Big| 
 \le C_2  c^n \dt 
 \le  C_2  C_1 \dt  .  \label{convergence-A-4-3}  
\end{eqnarray}   
Moreover, since  $q_1^{n+1} > 0$, $q_2^{n+1} > 0$, $q_3^{n+1} > 0$, it is straightforward to see that
	\begin{align}  
\Big| \frac{q_1^{n+1} e_A^n + q_2^{n+1} e_B^n - q_3^{n+1} e_C^n }{M^{n+1}} \Big| & \le   \frac{q_1^{n+1} + q_2^{n+1} + q_3^{n+1}}{M^{n+1}} \max ( |  e_a^n  | , | e_b^n | , | e_c^n | )    
	\nonumber 
	\\
&\le \max ( |  e_a^n  | , | e_b^n | , | e_c^n | )  .
	\label{convergence-A-5}   
	\end{align} 
A substitution of~\eqref{convergence-A-4-1}, \eqref{convergence-A-4-3} and \eqref{convergence-A-5} into \eqref{convergence-A-3-1} leads to 
	\begin{align} 
| e_R^{n+1} |  & \le   \frac{|q_0^{n+1} |}{M^{n+1}}  | e_c^n |  + \Big| \frac{q_1^{n+1} e_a^n + q_2^{n+1} e_b^n  - q_3^{n+1} e_c^n }{M^{n+1}}  \Big| + \frac{1}{M^{n+1}} | \tau_1^{n+1} |
	\nonumber  
	\\
&\le C_1  C_2 \dt  | e_c^n |  +  \max ( |  e_a^n  | , | e_b^n | , | e_c^n | ) +   C_1 \dt   | \tau_1^{n+1} |  . 
	\label{convergence-A-6}  
	\end{align} 
With the a-priori numerical error assumption at the previous time step \eqref{a priori-1}, we arrive at a rough error estimate for $e_R^{n+1}$:  
\begin{eqnarray} 
   | e_R^{n+1} |  &\le&  
    2 ( \dt^\frac12 + h) +   C_1 \dt   | \tau_1^{n+1} |   
  \le 2 ( \dt^\frac12 + h) +   C C_1 \dt^2   
  \le 3 \dt^\frac12 + 2 h ,  \label{convergence-A-7}  
\end{eqnarray}      
provided that $C_1 C_2 \dt \le 1$ and $C C_1 (\dt)^{3/2} < 1$. Here the local truncation error estimate $| \tau_1^{n+1} | \le C \dt$ has been used. 

The rough error estimate on $e_R^n$ enables us to refine the estimates on $q_i^{n+1}$, which is the key to obtain the error estimate of the desired order. As a result of this rough estimate, the following estimates can be derived: 
	\begin{eqnarray} 
&&  c_e^n + R_e^{n+1} \ge \epsilon_0 - C_0 \dt \ge \frac{\epsilon_0}{2} ,  \quad 
  \mbox{(since $| R_e ^{n+1} | \le C_0 \dt$)} ,  
  \label{convergence-A-8-1}    
\\
  && 
  | e_c^n | +  | e_R^{n+1} |  \le 4 \dt^\frac12 + 3 h \le \min\left(\frac{\epsilon_0}{4}, 1\right) ,  \quad 
  \mbox{using~\eqref{a priori-1}, \eqref{convergence-A-7})} ,  
  \label{convergence-A-8-2}   
\\
  && 
  c^n + R^{n+1}  \ge c_e^n + R_e^{n+1}  - ( | e_c^n | +  | e_R^{n+1} |  ) 
  \ge \frac{\epsilon_0}{4},   \label{convergence-A-8-3}      \\
  &&    c^n + R^{n+1}  \leq c_e^n + R_e^{n+1}  + ( | e_c^n | +  | e_R^{n+1} |  ) 
  \leq C_1 + 1,   \label{convergence-A-8-4}       
\end{eqnarray}
provided that $C_0 \dt \leq \frac{\epsilon_0}{2}$ and $4 \dt^\frac12 + 3 h \le \min(\frac{\epsilon_0}{4}, 1)$. The same estimate can be made for $a$ and $b$.
Then we obtain 
\begin{eqnarray} 
  0 < q_1^{n+1} = \frac{b^n - R^{n+1}}{c^n +R^{n+1}}  
  &\le& \frac{C_1 +1}{\frac{\epsilon_0}{4}} = 4 (C_1 +1) \epsilon_0^{-1} ,  
   \label{convergence-A-8-5} 
\\   
  0 < q_2^{n+1} = \frac{a_e^n - R_e^{n+1}}{c^n +R^{n+1}}    
  &\le& \frac{C_1}{\frac{\epsilon_0}{4}} = 4 C_1 \epsilon_0^{-1}  ,  
  \label{convergence-A-8-6}   
\\    
  0 < q_3^{n+1} = \frac{(a_e^n - R_e^{n+1}) (b_e^n - R_e^{n+1})}{(c^n + R^{n+1} )  
  (c_e^n +R_e^{n+1})} 
  &\le&  \frac{C_1^2}{\frac{\epsilon_0^2}{8}} = 8 C_1^2 \epsilon_0^{-2} , 
    \label{convergence-A-8-7}   
\end{eqnarray} 
so that the following uniform bound is available: 
\begin{equation} 
  0 < q_1^{n+1} + q_2^{n+1} + q_3^{n+1}  
  \le C_3 := ( 8 C_1 +4) \epsilon_0^{-1}   
  + 8 C_1^2 \epsilon_0^{-2}  .  \label{convergence-A-8-8}   
\end{equation}   

Consequently, we have the refined estimate
	\begin{eqnarray} 
\Big| \frac{q_1^{n+1} e_a^n + q_2^{n+1} e_b^n - q_3^{n+1} e_c^n }{M^{n+1}}  \Big| &\le&  \frac{1}{M^{n+1}}  ( q_1^{n+1} + q_2^{n+1} + q_3^{n+1} ) \max ( |  e_a^n  | , | e_b^n | , | e_c^n | )  
	\nonumber 
	\\
&\le& C_1 C_3 \dt \max ( |  e_a^n  | , | e_b^n | , | e_c^n | ) .  
   \label{convergence-A-8-9}  
	\end{eqnarray} 
Going back the earlier error estimate~\eqref{convergence-A-6}, we arrive at 
\begin{eqnarray} 
   | e_R^{n+1} |  &\le&   \frac{|q_4^{n+1} |}{M^{n+1}}  | e_c^n |   
   + \Big| \frac{q_1^{n+1} e_a^n + q_2^{n+1} e_b^n 
   - q_3^{n+1} e_c^n }{M^{n+1}}  \Big|    
   + \frac{1}{M^{n+1}} | \tau_1^{n+1} |    \nonumber  
\\
  &\le& 
  C_1  C_2 \dt  | e_c^n |   
  +   C_1 C_3 \dt   
     \max ( |  e_a^n  | , | e_b^n | , | e_c^n | )  
  +   C_1 \dt   | \tau_1^{n+1} |  \nonumber 
\\
  &\le&  2C_1  ( C_2 + C_3 ) \dt  \max ( |  e_a^n  | , | e_b^n | , | e_c^n | )  
  +   C_1 \dt   | \tau_1^{n+1} |  . 
   \label{convergence-A-9-1}  
\end{eqnarray} 
On the other hand, a difference between the numerical solution~\eqref{splitting-Field-A-3} and the constructed profile~\eqref{consistency-12-3} reveals that 
\begin{equation} 
     e_a^{n+1,*} =   e_a^n - e_R^{n+1} ,  \, \, \,  
     e_b^{n+1,*}  =  e_b^n - e_R^{n+1} , \, \, \,  
     e_c^{n+1,*}  =  e_c^n + e_R^{n+1} .   \label{convergence-A-9-2}     
\end{equation} 
Then we arrive at the following error estimate in the first stage: 
\begin{eqnarray}  
  && 
   | e_a^{n+1,*} | \le | e_a^n | + | e_R^{n+1} | 
   \le ( 1+ C_4 \dt )  \max ( |  e_a^n  | , | e_b^n | , | e_c^n | )  
  +   C_1 \dt   | \tau_1^{n+1} |   ,    \label{convergence-A-9-3}  
\\
  && 
   | e_b^{n+1,*} | \le | e_b^n | + | e_R^{n+1} | 
   \le ( 1+ C_4 \dt )  \max ( |  e_a^n  | , | e_b^n | , | e_c^n | )  
  +   C_1 \dt   | \tau_1^{n+1} |   ,    \label{convergence-A-9-4}   
\\
  && 
   | e_c^{n+1,*} | \le | e_c^n | + | e_R^{n+1} | 
   \le ( 1+ C_4 \dt )  \max ( |  e_a^n  | , | e_b^n | , | e_c^n | )  
  +   C_1 \dt   | \tau_1^{n+1} |   ,    \label{convergence-A-9-5}        
\end{eqnarray}  
with $C_4 :=  ( 1 + M_0 ) C_1  ( C_2 + C_3 )$. Since the error estimate~\eqref{convergence-A-9-3}-\eqref{convergence-A-9-5} is valid at a point-wise level, the following conclusion is made: 
\begin{equation}  
  \begin{aligned} 
    & 
   \| e_a^{n+1,*} \|_\infty , \,  \| e_b^{n+1,*} \|_\infty , \, \| e_c^{n+1,*} \|_\infty  
\\
   \le & ( 1+ C_4 \dt )  \max ( \|  e_a^n  \|_\infty , \| e_b^n \|_\infty , 
   \| e_c^n \|_\infty )  +   C_1 \dt   \| \tau_1^{n+1} \|_\infty  .  
  \end{aligned}    
  \label{convergence-A-9-6}   
\end{equation} 

\begin{remark} 
In the rough error estimate~\eqref{convergence-A-7}, we see that the accuracy order is lower than the desired accuracy order. Therefore, such a rough estimate could not be used for a global induction analysis. Instead, the purpose of such an estimate is to establish a uniform $\| \cdot \|_\infty$ bound, so that a discrete separation property becomes available for the numerical solution, as well as its maximum values. With such a property established for the numerical solution, the refined error analysis yields much sharper estimate as in~\eqref{convergence-A-9-1}. A combination of a rough error estimate and a refined error estimate has been successfully applied to certain nonlinear PDEs with singular terms, such as the Poisson-Nernst-Planck system~\cite{liu2020positivity}, the  porous medium equation in the energetic variational formulation~\cite{duan20a}. Here we show that such a technique works for the highly nonlinear reaction trajectory equation (\ref{Equation_R})
\end{remark} 

\subsection{Error estimate in the second stage} 
Now we proceed into the error estimate for the second part. Subtracting the implicit Euler scheme~\eqref{splitting-Field-B-1} from the consistency estimate~\eqref{consistency-12-4} yields
\begin{eqnarray} 
  & 
  \frac{e_a^{n+1} - e_a^{n+1,*}}{\dt} = 
  \nabla_h \cdot (  D_a \nabla_h  e_a^{n+1} )  
   + \tau_{2,a}^{n+1} ,  
  \label{error equation-B-1}    
\\
  &
    \frac{e_b^{n+1} - e_b^{n+1,*}}{\dt} = 
    \nabla_h \cdot (  D_b \nabla_h  e_b^{n+1}  )  
    + \tau_{2,b}^{n+1} ,  
  \label{error equation-B-2}  
\\
  &
 \frac{e_c^{n+1} - e_c^{n+1,*}}{\dt} = 
  \nabla_h \cdot (  D_c \nabla_h  e_c^{n+1}  )   
   + \tau_{2,c}^{n+1} ,  
  \label{error equation-B-3} 
 \end{eqnarray} 
where the local truncation errors $\tau_{2,a}^{n+1}$, $\tau_{2,b}^{n+1}$ and $\tau_{2,c}^{n+1}$ satisfy $| \tau_{2, a}^{n+1} |, | \tau_{2,b}^{n+1} | , | \tau_{2,c}^{n+1} | \le C (\dt + h^2)$, at a point-wise level. 

Due to the maximum principle for the discrete elliptic operator in the finite difference setting \cite{isaacson2012analysis}, we have 
	\begin{eqnarray}  
  && 
   \| e_a^{n+1} \|_\infty  \le \| e_a^{n+1,*} \|_\infty + \dt \| \tau_{2,a}^{n+1} \|_\infty , 
   \label{convergence-B-1-1}  
\\
  && 
   \| e_b^{n+1} \|_\infty  \le \| e_b^{n+1,*} \|_\infty + \dt \| \tau_{2,b}^{n+1} \|_\infty , 
   \label{convergence-B-1-2}  
\\
  && 
   \| e_c^{n+1} \|_\infty  \le \| e_c^{n+1,*} \|_\infty + \dt \| \tau_{2,c}^{n+1} \|_\infty .
   \label{convergence-B-1-3}        
	\end{eqnarray}    
Indeed, for equation (\ref{error equation-B-1}) with periodic boundary condition, if $e_a^{n+1}$ takes a maximum value at $(i,j,k)$, we see that 
	\begin{equation}
\nabla_h \cdot (  D_a  \nabla_h  e_a^{n+1} )_{i,j,k}  \le 0,
	\end{equation}
by looking at the values of $e_a^{n+1}$ in a neighborhood of $(i, j, k)$, provided that $ D_a ({\bf x}) $ is point-wise non-negative. Therefore, the following inequality is valid 
	\begin{equation}
(e_a^{n+1})_{i,j,k} \le ( e_a^{n+1,*})_{i,j,k} + \dt ( \tau_{2,a}^{n+1} )_{i,j,k},
	\end{equation}
which in turn implies that
\begin{equation}\label{dmp_1}
\max_{i,j,k} e_a^{n+1} \le \max_{i,j,k} e_a^{n+1,*}  + \dt \max_{i,j,k} \tau_{2,a}^{n+1} . 
\end{equation}
Similarly, we can prove that
\begin{equation}\label{dmp_2}
\min_{i,j,k} e_a^{n+1} \ge \min_{i,j,k} e_a^{n+1,*}  - \dt \max_{i,j,k} | \tau_{2,a}^{n+1} |.
\end{equation}
Combining (\ref{dmp_1}) and (\ref{dmp_2}), we obtain (\ref{convergence-B-1-1}). Inequalities (\ref{convergence-B-1-2}) and (\ref{convergence-B-1-3}) can be proved in the same manner.

\subsection{Convergence estimate for the full operator splitting system} 

A combination of~\eqref{convergence-A-9-6} and \eqref{convergence-B-1-1}-\eqref{convergence-B-1-3} reveals that 
\begin{eqnarray}  
   \max ( \| e_a^{n+1} \|_\infty , \| e_b^{n+1} \|_\infty , \| e_c^{n+1} \|_\infty ) 
   &\le& ( 1+ C_4 \dt )  \max ( \|  e_a^n  \|_\infty , \| e_b^n \|_\infty , 
   \| e_c^n \|_\infty )   \nonumber 
\\
  && 
   +   ( 1+ C_1 ) \dt  (  \| \tau_1^{n+1} \|_\infty +  \| \tau_2^{n+1} \|_\infty ) .    
    \label{convergence estimate-1}   
\end{eqnarray}  
Therefore, an application of a discrete Gronwall inequality leads to the desired convergence estimate
	\begin{equation}
\max ( \| e_a^{n+1} \|_\infty , \| e_b^{n+1} \|_\infty , \| e_c^{n+1} \|_\infty ) 
   \le C ( \dt + h^2 ) , \label{convergence estimate-2}
	\end{equation}  
based on the truncation error estimates $\| \tau_1^{n+1} \|_\infty \le C \dt$,   $\| \tau_2^{n+1} \|_\infty \le C ( \dt + h^2 )$. 

With the $\| \cdot \|_\infty$ error estimate~\eqref{convergence estimate-2} at hand, the a-priori assumption in~\eqref{a priori-1} is satisfied at the next time step $t^{n+1}$:  
	\begin{equation} 
  \| e^{n+1} \|_\infty  \le C ( \dt + h^2 ) \le \dt^\frac12 + h , 
	\label{a priori-3}  
	\end{equation} 
provided $\dt$ and $h$ are sufficiently small. As a result, an induction analysis could be applied. This finishes the proof of Theorem~\ref{thm:convergence}.

	\begin{remark} 
There have been many existing works of operator splitting numerical approximation to nonlinear PDEs, such as \cite{descombes2001convergence, descombes2004operator, zhao2011operator} for reaction-diffusion systems, \cite{Bao2002, besse2002order, Lubich2008, Shen2013, Thalhammer2012} for the nonlinear Schr\"odinger equation, \cite{Badia2013} for the incompressible magnetohydrodynamics system, \cite{BCF2013} for the delay equation, \cite{EO2014a} for the nonlinear evolution equation, \cite{EO2014b} for the Vlasov-type equation, \cite{KV2015} for a generalized Leland's mode, \cite{zhangC18a, zhangC17a} for the ``Good" Boussinesq equation, \cite{Lee2015} for the Allen-Cahn equation,  \cite{Li2017} for the molecular beamer epitaxy (MBE) equation, \cite{Zhao2014} for nonlinear solvation problem, etc. A few convergence estimates have also been reported for gradient flow with polynomial energy potential, such as~\cite{Li2017, zhangC17a}. 
 The convergence result stated in this article provides a theoretical convergence analysis for an operator splitting scheme for an energy variational formulation with singular energy potential involved. 
	\end{remark}

\section{Numerical test}
\begin{figure}[!h]
   \begin{overpic}[width = 0.48\linewidth]{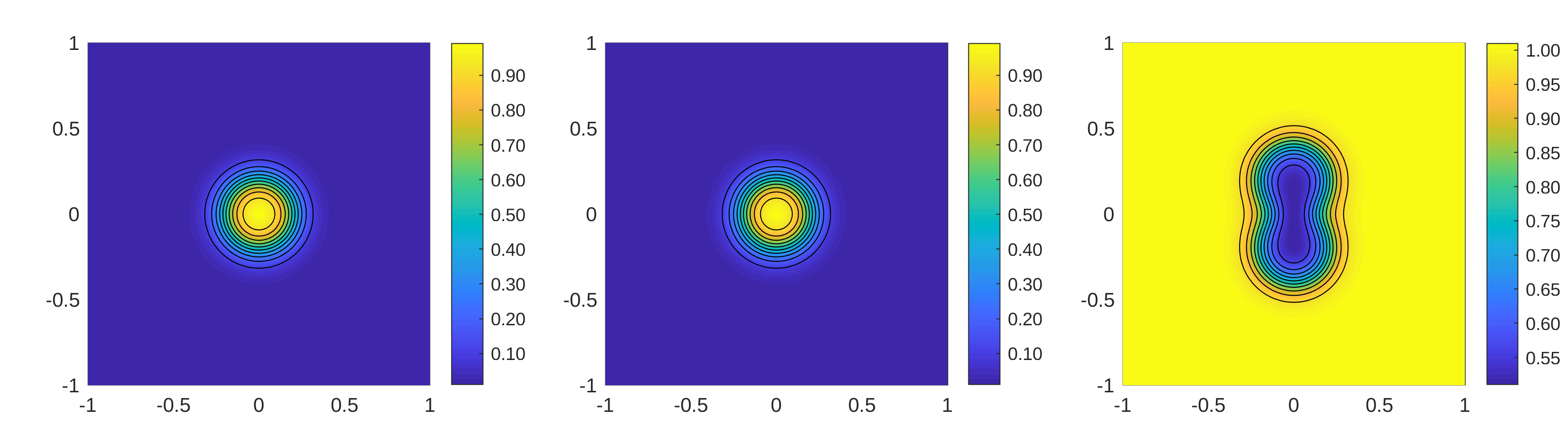}
    \put(-3, 22){ \scriptsize (a)}
   \end{overpic}
   \hfill
   \begin{overpic}[width = 0.48\linewidth]{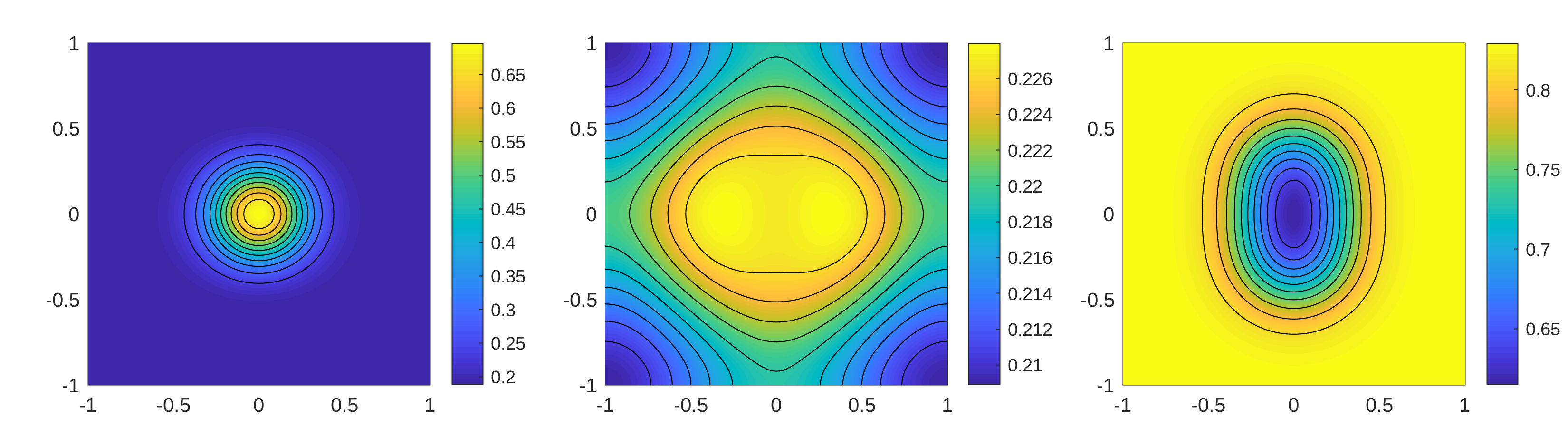}
   \put(-3, 22){ \scriptsize (b)}
  \end{overpic}

  \vspace{0.5 em}
  \begin{overpic}[width = 0.48\linewidth]{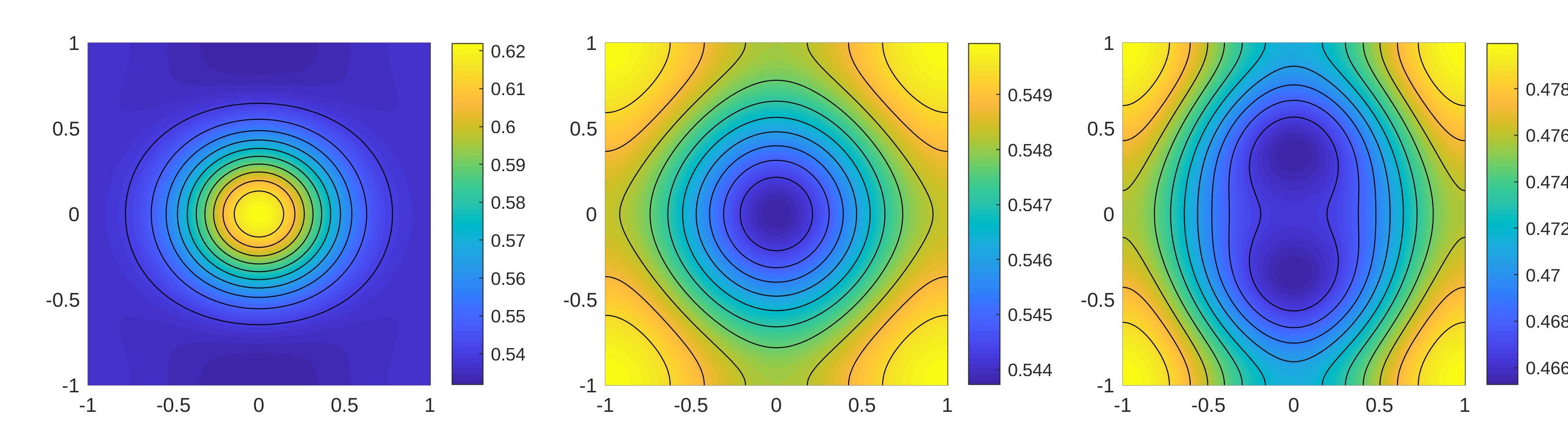}
   \put(-3, 22){ \scriptsize (c)}
  \end{overpic}
   \hfill
   \begin{overpic}[width = 0.48\linewidth]{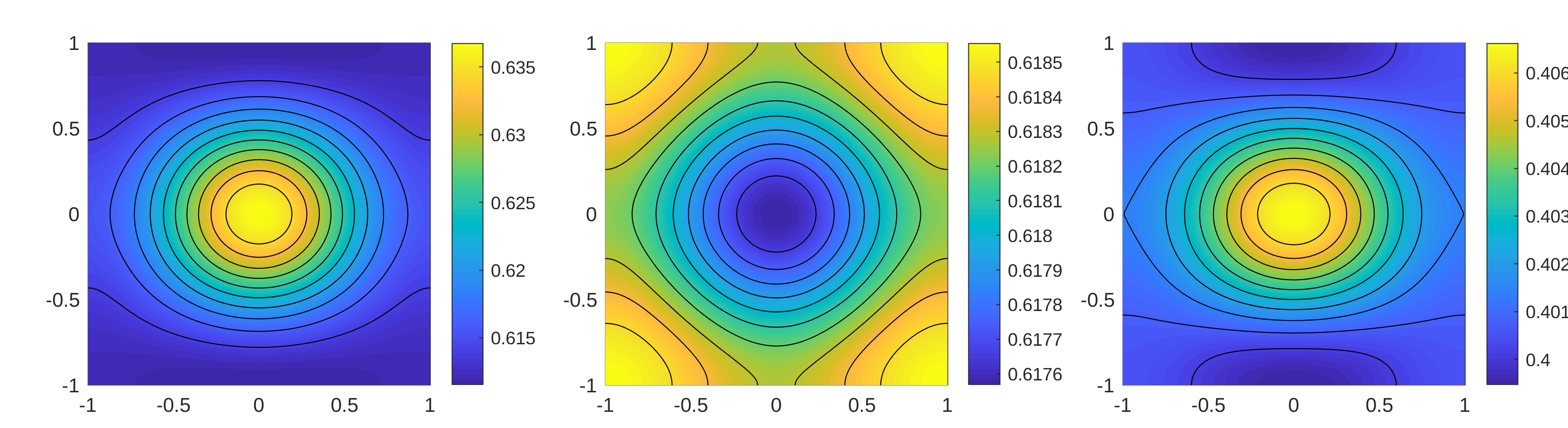}
   \put(-3, 22){ \scriptsize (d)}
  \end{overpic}
  \caption{(a)-(d): Numerical results for the reaction-diffusion system  \eqref{Reaction-D} with $D_a \equiv 0.05$, $D_b \equiv 1$, $D_c \equiv 0.1$ and the initial condition (\ref{initial}): (a) t = 0, (b) t = 0.2, (c) t = 1 and (d) t = 2.}\label{RD_1}
	\end{figure}

In this section, we present a 2D numerical example for equation \eqref{Reaction-D}. The computational domain is taken as $\Omega = (-1, 1)^2$, and the initial condition is set as
\begin{equation}\label{initial}
\begin{cases}
& a_0(x, y) = \frac{1}{2}(- \tanh( \frac{\sqrt{x^2 + y^2} - 0.2}{0.1}) + 1) + 0.01; \\
& b_0(x, y) = \frac{1}{2}( \tanh( \frac{\sqrt{x^2 + y^2} - 0.2}{0.1}) + 1) + 0.01; \\
& c_0(x, y) = \frac{1}{4}   \tanh( \frac{\sqrt{x^2 + (y - 0.2)^2} - 0.2}{0.1} + 1) + \tfrac{1}{4}\tanh(  \tfrac{\sqrt{x^2 + (y + 0.2)^2} - 0.2}{0.1} + 1) + 0.01 .
\end{cases} 
\end{equation}
The diffusion coefficients are given by $D_a \equiv 0.05$, $D_b \equiv 1$ and $D_c \equiv 0.1$. The initial condition and numerical solutions at different time instants are displayed in Figure~\ref{RD_1},

In addition, we look at the numerical error at $T = 0.2$, before the system reaches the constant steady state. Since the analytical solution is not available, we use the numerical solution with $h = 1/200$ and $\Delta t = 1/1600$ as the reference solution in the accuracy test for the temporal numerical errors. Moreover, we fix the spatial resolution as $h=\frac{1}{200}$ for the temporal accuracy test, so that the spatial numerical error is negligible. Table~\ref{table2} displays the $\| \cdot \|_{\infty}$ numerical errors at $T = 0.2$ with a sequence of time step sizes: $\dt=\frac{1}{25}$, $\frac{1}{50}$, $\frac{1}{100}$, $\frac{1}{200}$ and $\frac{1}{400}$. This result has which indicates a clear first order accuracy in time.
\begin{table}[!h]
  \begin{center}
    \begin{tabular}{c | c | c |c | c | c | c | c }
     \hline \hline
     $\dt$ &  h & $\| e_a \|_{\infty}$ & Order &  $\| e_b \|_{\infty}$ & Order &  $\| e_c \|_{\infty}$ & Order \\
     \hline
     1/25   &  1/200  & 9.5498e-3  &         & 1.2498e-2  &         & 7.1119e-3  &       \\
     1/50   &  1/200  & 4.8519e-3  & 0.9769  & 5.8081e-3  & 1.1056  & 3.5450e-3  & 1.0044 \\
     1/100  &  1/200  & 2.3840e-3  & 1.0252  & 2.7387e-3  & 1.0846  & 1.7314e-3  & 1.0338 \\
     1/200  &  1/200  & 1.1208e-3  & 1.0889  & 1.2629e-3  & 1.1168  & 8.1173e-4  & 1.0929 \\
     1/400  &  1/200  & 4.8213-4   & 1.2170  & 5.3817e-4  & 1.2306  & 3.4862e-4  & 1.2193 \\
     \hline  \hline
    \end{tabular}

\end{center}
\caption{Numerical errors, order of accuracy for numerical simulations of (\ref{Reaction-D}) with $D_a = 0.05$, $D_b = 1$, $D_c = 0.1$ and initial condition (\ref{initial})  at $T = 0.2$. The numerical solution with $h = 1/200$ and $\dt = 1/1600$ is taken as the reference solution.}\label{table2}
\end{table}

To test the spatial accuracy of the operator scheme for this example, we perform the computations on a sequence of mesh resolutions: $h=\frac{1}{20}, \frac{1}{30}, \frac{1}{40}, \frac{1}{50}$, $\frac{1}{60}$, and the time step size is set as $\Delta t=h^2$ to eliminate the affect of temporal errors.
Since an analytical form of the exact solution is not available, we compute the $\ell^{\infty}$ differences between numerical solutions with consecutive spatial resolutions, $h_{j-1}$, $h_j$ and $h_{j+1}$, in the Cauchy convergence test. 
Since we expect the numerical scheme preserves a second order spatial accuracy, we can compute the following quantity
$$
   \frac{ \ln \Big(  \frac{1}{A^*} \cdot 
   \frac{\| u_{h_{j-1}} - u_{h_j} \|_\infty }{ \| u_{h_j} - u_{h_{j+1}} \|_\infty} \Big) } 
   {\ln  \frac{h_{j-1}}{h_j} } ,  \quad A^* =  \frac{ 1 - \frac{h_j^2}{h_{j-1}^2} }{1 - \frac{h_{j+1}^2}{h_j^2} } ,  
   \quad \mbox{for} \, \, \, h_{j-1} > h_j > h_{j+1} , 
$$
to check the convergence order \cite{liu2020positivity}. As demonstrated in Table~\ref{t2:convergence}, 
an almost perfect second order spatial convergence rate for the proposed operator splitting scheme is observed. 

\begin{table}[ht]
  \begin{center}
  \begin{tabular}{c |c |c |c |c |c |c }
  \hline  \hline
   --- & $\psi = a$ & Order & $\psi = b$& Order & $\psi = c$& Order \\
   \hline
  $\| \psi_{h_1} - \psi_{h_2} \|_{\infty}$  & 2.0358e-3 & -       & 4.1584e-4  & -      & 7.6602e-4 & - \\
  $\| \psi_{h_2} - \psi_{h_3} \|_{\infty}$  & 7.1819e-4 & 1.9805  & 1.4459e-4  & 2.0162  & 2.6167e-4 & 2.0599  \\ 
  $\| \psi_{h_3} - \psi_{h_4} \|_{\infty}$  & 3.3291e-4 & 1.9949  & 6.6751e-5  & 2.0090  & 1.2073e-4 & 2.0111   \\
  $\| \psi_{h_3} - \psi_{h_4} \|_{\infty}$  & 1.8086e-4 & 1.9995  & 3.6211e-5 & 2.0060  & 6.5512e-5 & 2.0048    \\
   \hline  \hline
  \end{tabular}
  \caption{The $\ell^\infty$ differences and convergence order for the numerical solutions of $a$, $b$, and $c$ at $T=0.2$. Various mesh resolutions are used: $h_1=\frac{1}{20}$, $h_2=\frac{1}{30}$, $h_3=\frac{1}{40}$, $h_4=\frac{1}{50}$, $h_5=\frac{1}{60}$, and the time step size is taken as $\Delta t=h^2$.}
  \label{t2:convergence}
  \end{center}
  \end{table}

\section{Concluding remarks} \label{sec: conclusion} 


  A detailed convergence analysis and error estimate have been presented for the variational operator splitting scheme for the reaction-diffusion system~\eqref{Reaction-D}, which satisfies the detailed balance condition. The operator splitting scheme is based on an energetic variational formulation, in which the equation of the reaction trajectory $R$ is introduced in the reaction stage, and both the reaction and diffusion stages dissipate the same discrete free energy \cite{liu2020structure}.
   To overcome a well-known difficulty associated with the implicit treatment of the highly nonlinear and singular nature of the logarithmic terms, we make use of the convex nature of these nonlinear terms  
   A combination of rough error estimate and refined error estimate leads to a desired bound of the numerical error at the reaction stage, in the discrete $\| \cdot \|_\infty$ norm. In addition, a discrete maximum principle yields the evolution bound of the numerical error function at the diffusion stage. As a result, a combination of the numerical error analysis at different stages and the consistency estimate for the operator splitting yields the desired convergence estimate for the full reaction-diffusion system in the discrete $\| \cdot \|_\infty$ norm, provided that the exact solution are sufficiently smooth, and $\dt$ and $h$ are sufficiently small. It is straightforward to extend the analysis to other reaction-diffusion systems with detailed balance condition \cite{liu2020structure}.
   
\section*{Acknowledgement} 
This work is partially supported by the National Science Foundation (USA) grants NSF DMS-1759536, NSF DMS-1950868 (C. Liu, Y. Wang), NSF DMS-2012669 (C. Wang), and NSF DMS-1719854, DMS-2012634 (S.~Wise).  Y. Wang would also like to thank Department of Applied Mathematics at Illinois Institute of Technology for their generous support and for a stimulating environment.

\bibliographystyle{siam}
\bibliography{KCR}

\end{document}